\documentstyle[12pt,amssymbols]{article}  
\newcommand{\bq}{\begin{equation}}  
\newcommand{\eq}{\end{equation}}  
\newcommand{\bqa}{\begin{eqnarray}}  
\newcommand{\eqa}{\end{eqnarray}}  
\newcommand{\ra}{\rightarrow}  
  
\def\L{\Lambda}

\def\half{{1 \over 2}}  
\def\s{\sigma}

\def\ep{\epsilon}

\def\ov{\over}

\def\ed{\end{document}}  
  
\def\ra{\rightarrow}  
\def\al{\alpha}  
\def\2pi{1\over 2\pi i}  

\def\~{\tilde}
\def\newline{\hfil\break}

\def\ra{\rightarrow}

\def\sq2{\sqrt{2}}  
\def\sqk2{\sqrt{2(k+2}}  
\def\sqk{\sqrt{k}}

\def\be{\begin{equation}}  
\def\ee{\end{equation}}  
\def\br{\begin{array}}  
\def\er{\end{array}}  
\def\bea{\begin{eqnarray}}  
\def\eea{\end{eqnarray}}  
\def\ba{\begin{equation}\begin{array}}  
\def\ea{\end{array}\end{equation}}  
\def\bac{\begin{equation}\begin{array}{rll}}

\def\al{\alpha}

\newcommand{\uq}{U_q (\widehat{sl(2)})}

\def\Z{{\Bbb Z}}  
  
\def\N{{\Bbb N}}  
\def\L{{\Bbb L}}
\def\R{{\Bbb R}}

\def\pl{\prod\limits}

\def\ep{\epsilon}

\begin{document}  
\rightline{ITP-SB-96-41}  
\rightline{August, 1996}  
\vbox{\vspace{-10mm}}  
\vspace{1.0truecm}  
\begin{center}  
{\LARGE \bf   
Exact  $n$-spinon dynamic 
structure factor
of the Heisenberg model
 }\\[8mm]  
{\large A.H. Bougourzi
}\\  
[6mm]{\it 
Institute of Theoretical Physics\\  
SUNY at Stony Brook\\  
Stony Brook, NY 11794
}\\[20mm]  

\end{center}  
\vspace{1.0truecm}  
\begin{abstract}  
Using the spinon picture we derive an  integral 
representation for the exact 
$n$-spinon dynamic structure factor
of the spin-1/2  Heisenberg model. 
\end{abstract}  
\newpage  
\section{Introduction}  

Since the pionnering work of Bethe \cite{Bet31}, 
the Heisenberg model in the anti-ferromagnetic (AF) regime 
has received considerable attention \cite{Hul38}-\cite{
Fleal96}. 
It is also of great experimental interest where 
the dynamic structure factor (DSF) plays a central role in 
magnetic neutron scattering experiments \cite{Ten95,Tenal95}.
We extend the recent derivation  of the 2-spinon DSF of the 
$s=1/2$ AF  
Heisenberg model  
to the exact $n$-spinon one \cite{Boual96}. For an extensive
analysis of this 2-spinon DSF see Refs. \cite{Karal96,Fledd96}. 
Though more complicated, the result
of this paper can be thought of as the analogue of 
that of the $s=1/2$ $XX$ model at zero temperature, and which
was first derived by Niemeijer \cite{Nie67} and later
re-derived in a very simple manner by
 M\"uller at al \cite{Mulal81}. 
Our work is based on the  spinon picture of the eigenspace of
the Heisenberg model in the AF regime, 
which was rigorously 
established in Refs.
\cite{FaTa79,Babal83}. 
We use the form factors that  
have been derived in Ref. \cite{JiMi94} from the isotropic 
limits of those of the XXZ model.

This letter is organized as follows:
first, we briefly review the 
spinon picture of the eigenspace of the $s=1/2$ XXZ Hamiltonian.
Then, we define one component of the DSF for this model. 
By taking the isotropic limit,  we derive an exact integral 
representation for 
all non-vanishing components of the 
$n$-spinon 
DSF. The summation over all these
$n$-spinon contributions leads to the total DSF of the AF
Heisenberg 
model.

\section {The spinon picture in the XXZ model}

The Hamiltonian of the anisotropic (XXZ) Heisenberg model 
is defined by
\be { H_{XXZ}}=
-\half \sum_{n=-\infty}^{\infty} (\s_n^x \s_{n+1}^x 
+\s_n^y \s_{n+1}^y + \Delta\s_n^z \s_{n+1}^z) 
\label{hamiltonian}, 
\ee 
where $\Delta=(q+q^{-1})/2$ is the anisotropy parameter. 
Here $\sigma_n^{x,y,z}$ 
are the usual
Pauli matrices acting at the $n^{\rm th}$ position of the 
formal 
infinite tensor product 
\be W= \cdots V \otimes V \otimes V  \cdots 
\label{infprod},\eq 
where $V$ is the two-dimensional representation of $U_q(sl(2))$ 
quantum group. 
We consider the model in 
the anti-ferromagnetic regime $\Delta <-1$, i.e., 
$-1<q<0$. Later we take the isotropic limit $q\ra -1$ in 
a special manner.
This model is symmetric 
under the quantum group $\uq$, and therefore its 
 eigenspace
is identified with the following level 0 $\uq$ 
module \cite{Daval93}:
\be
{\cal F}=\sum_{i,j}V(\Lambda_i)\otimes V(\Lambda_j)^*. 
\ee
Here $\Lambda_i$ and $V(\Lambda_i); i=0,1$ are  level 
1 $\uq$-highest
weights and $\uq$-highest weight modules, respectively, whereas
$V(\Lambda_i)^*$ are  dual modules defined from  
$V(\Lambda_i)$  
through the antipode.
The module $V(\Lambda_i)$ is identified 
with the
subspace of the following formal semi-infinite tensor 
product of $V$'s:
\be X= \cdots V \otimes V \otimes V,\eq 
and which consists  of all linear combinations of  
spin configurations with  fixed boundary conditions such that 
the eigenvalues of 
$\sigma^z_n$ are $(-1)^{i+n}$ in the limit $
 n\ra -\infty$. Similarly, the module $V(\Lambda_i)^*$ is 
associated with the right semi-infinite tensor product of the
$V$'s. The  
spinon picture 
of the eigenspace of 
this Hamiltonian is given in terms of vertex operators 
(of type II according to the terminology of Ref. \cite{Daval93}) 
which
act as intertwiners of $\uq$ modules, and which create 
the set of eigenstates (i.e., spinons) $\{|\xi_1,\cdots 
\xi_n>_{\ep_1,\cdots \ep_n;i}, n\geq 0\}$. Here $i=0,1$ fixes
the boundary conditions, $\xi_j$ denotes  a spectral parameter
living on the unit circle,
and $\ep_j=\pm 1$ is twice the $z$-component of   
the spin of a
spinon. The actions on ${\cal F}$ of $H_{XXZ}$ and the 
translation operator 
$T$, which shifts the spin  chain by one site,  
are given by 
\bac T|\xi_1,\cdots,\xi_n>_i &=&\pl_{i=1}^n\tau(\xi_i)^{-1}
|\xi_1,\cdots,\xi_n>_{1-i},\quad 
T|0>_i =|0>_{1-i}, \\
H_{XXZ}|\xi_1,\cdots,\xi_n>_i 
&=&\sum_{i=1}^n e(\xi_i)|\xi_1,\cdots,\xi_n>_i, 
\label{states}\ea  
where
\bac \tau(\xi)&=& \xi^{-1} 
{\theta_{q^4} (q \xi^2) \ov \theta_{q^4} (q 
\xi^{-2})}=e^{-i p(\alpha)},\quad p(\alpha)={\rm am}({2K\over \pi}
\alpha)+{\pi/2}, \\
e(\alpha) &=&{ 1-q^2 \ov 2 q} \xi {d \ov d \xi} 
\log \tau(\xi)= {2K\over \pi} \sinh({\pi K^\prime\over K})
{\rm dn}({2K\over \pi}\alpha). \label{enmom}
\ea
Here, $e(\alpha)$ and  $p(\alpha)$ are the energy and the 
momentum of
the spinon respectively,  ${\rm am}(x)$ and ${\rm dn}(x)$ 
are the usual 
elliptic amplitude and
delta functions, with nome $-q$ and 
 complete elliptic integrals
 $K$ and $K^\prime$, and 
\bac
q&=&-\exp(-\pi K^\prime/K),\\
\xi&=&ie^{i\alpha},\\ 
\theta_x(y)&=&(x;x)_{\infty} (y;x)_{\infty} 
(x y^{-1};x)_{\infty},\\
(y;x)_{\infty}&=&\prod_{n=0}^{\infty} (1-y x^n).
\ea
Therefore, $\sigma^{x,y,z}(t, n)$ at time $t$ and position
$n$ are related to $\sigma^{x,y,z}(0,0)$ at time 0 and position
0 through:
\be \sigma^{x,y,z} (t,n)=\exp(i tH_{XXZ} ) T^{-n} 
\sigma^{x,y,z} 
(0,0) T^{n} 
\exp(-i tH_{XXZ} ). 
\ee
Moreover,  the completeness relation in ${\cal F}$ reads as 
\cite{JiMi94}: 
\be {\Bbb I}=\sum_{i=0,1}\sum_{n \geq 0} \sum_{\ep_1,\cdots,
\ep_n=\pm 1}
{1 \ov {n !}} \oint  {d\xi_1\over 2\pi i \xi_1} \cdots 
 {d\xi_n\over 2\pi i \xi_n}
 |\xi_n,\cdots,\xi_1>_{{\ep_n,\cdots,\ep_1};i}\;
{_{i;{\ep_1,\cdots,\ep_n}}{<\xi_1,\cdots,\xi_n|}}. 
\ee

\section{An integral representation of the $n$-spinon 
DSF of the 
Heisenberg model}

Let us first recall the definition of one of the components of
 the DSF in the case of the
XXZ
model in the sector $i$ as:
\be
S^{i,+-}(w, k)=
\int_{-\infty}^{\infty} dt \sum_{m\in\Z}
e^{i(wt+km)} {_i}< 0|\sigma^+(t, m)\sigma^-(0,0)|0>_i=\sum_{n\in 
2\N}S_n^{i,+-}(w,k),
\ee
where $w$ and $k$ are the energy and momentum transfer
respectively. 
It is sufficient to consider just the above component because 
in the isotopic limit, which is our main purpose later, 
all non-vanishing components are equal and proportional to 
$S^{i,+-}(w,k)$.
Due to the conservation of the total spin, 
only pairs of spinons 
contribute to the DSF. 
Using the completeness relation, 
the $n$-spinon contribution is given by
\bac
S_n^{i,+-}(w,k) &=& {2\pi\over n!} 
\sum_{m\in \Z} \sum_{\ep_1,\dots,\ep_n} 
\oint {d\xi_1\over 2\pi i \xi_1}
\dots \oint {d\xi_n\over 2\pi i \xi_n}
\exp\left(im(k+\sum_{j=1}^n p(\xi_j))\right)
 \delta (w-\sum_{j=1}^n e(\xi_j))\\
&&\times {_{i+m}<0|}\sigma^+ 
(0,0)|\xi_n,\dots,\xi_1>_{\ep_n,\dots,\ep_1;i+m}\:
{_{i;\ep_1,\dots,\ep_n}<\xi_1,\dots,\xi_n|}
\sigma^-(0,0)|0>_i,
\ea
which, in turn, can be re-written as
\bac
&& S_n^{i,+-}(w,k)={2\pi\over n!}
\sum_{\ep_1,\dots,\ep_n} 
 \oint {d \xi_1\over 2\pi i \xi_1}\dots\oint
{ d \xi_n\over 2 \pi i \xi_n} 
  \sum_{m\in \Z} \exp\left( 2 i
m(k+\sum_{j=1}^n p(\xi_j))\right)
\\
&&\times \delta (w-\sum_{j=1}^n e(\xi_j))\Bigl(
{_{i}<0|}\sigma^+(0,0)|\xi_n,\dots,\xi_1>_{\ep_n,\dots,\ep_1;i}
\:
{_{i;\ep_1,\dots,\ep_n}<\xi_1,\dots,\xi_n|}
\sigma^-{(0,0)}|0>_i \Bigr.\\
&& +
\Bigl. \exp\left( i(k+\sum_{j=1}^n p(\xi_j)\right)
{_{1-i}<0|}\sigma^+(0,0)|\xi_n,\dots,\xi_1>_{\ep_n,\dots,
\ep_1;1-i}\:
{_{i;\ep_1,\dots,\ep_n}<\xi_1,\dots,\xi_n|}
\sigma^-{(0,0)}|0>_i \Bigr).
\label{corrn}\ea
The non-vanishing form factors have
been computed in \cite{JiMi94}, and satisfy the
following relations:
\bac {_i{<}}0|\sigma^-(0,0)|\xi_n,\dots,\xi_1>_{\ep_n,\dots,
\ep_1;i} &=& 
{_{1-i}{<0|}}\sigma^+(0,0)|\xi_n,\dots,\xi_1>_{
-\ep_n,\dots,
-\ep_1;1-i}\\
&=&{_{i}{<0|}}\sigma^+(0,0)|-q\xi_1^{-1},\dots,-q\xi_n^{-1}>_{
-\ep_1,\dots,
-\ep_n;i}, \\
{_{i;\ep_1,\dots,
\ep_n}{<}} \xi_1,\dots,\xi_n|\sigma^-(0,0)|0>_i&=&
{_i{<0|}}\sigma^-(0,0)|-q \xi_1,\dots,-q \xi_n>_{
-\ep_1,\dots,
\ep_n;i}.
\label{FF}\ea

Their isotropic limits as $q\ra -1$ are performed
by first  making the following redefinitions:
\bac
\xi&=&ie^{{\epsilon\beta\over i\pi}},\\
q&=&-e^{-\epsilon},\quad\epsilon\ra 0^+,
\ea
with $\beta$, the appropriate
spectral parameter  for the
Heisenberg model, being real.

Then, one finds the following exact isotropic 
limits \cite{JiMi94}:\footnote{We do not find the same
overall coefficient in this limit as that of Ref. \cite{JiMi94}}
$^{,}$\footnote{We thank Karbach and M\"uller for confirming to
us the relation $|A_{-}(i\pi/2)|^2
|A_{+}(i\pi/2)|^2=1/2$, which is used here. 
See \cite{JiMi94} for the definition of 
$A_{+}(x)$} 
\bac
&&|{_{i}{<0|}}\sigma^+(0,0)|\xi_n,\dots,\xi_1>_{\ep_n,\dots,
\ep_1;i}|^2
{d\xi_{1}\over 2\pi i  \xi_{1}}\dots 
{d\xi_{n}\over  2\pi i \xi_{n}}\ra\\
&&{
2^{3n(1-n)/2}  |g_{\ep_1,\dots,\ep_n}
(\beta_1,\dots,\beta_n)|^2 
\prod_{1\leq j<j^\prime\leq
n}|A_{-}(\beta_{j^\prime}-
\beta_j)|^2\over \pi^{(7n^2-14n+8)/4}
\Gamma(1/4)^{n(n-2)}
|A_{-}(i\pi/2)|^{n(n-2)}
\prod_{j=1}^n \cosh(\beta_j)} d\beta_1\dots d\beta_n,\\
&&\lim({_{1-i}{<0|}}\sigma^+(0,0)|\xi_n,\dots,\xi_1>_{\ep_n,
\dots,\ep_1;1-i})=-
\lim({_{i}{<0|}}\sigma^+(0,0)|\xi_n,\dots,\xi_1>_{\ep_n,\dots,
\ep_1;i}),\\
&&p(\xi_j)\ra p(\beta_j)=\cot^{-1}(\sinh(\beta_j)),
\quad -\pi\leq p(\beta_j)
\leq 0,\\
&&e(\xi_j)\ra e(\beta_j)={\pi\over \cosh(\beta_j)}=
-\pi \sin(p(\beta_j)),\quad 1\leq j\leq n,
\ea
where
\bac
|A_{-}(x+iy)|^2&=&
\exp\left( -\int_{0}^{\infty}dt {(\cosh(2 t(1-{y\over\pi}))
\cos({2 t x
\over \pi})-1)
\exp(t)\over t \sinh(2 t)\cosh(t)}\right),\quad 
x,y\in\R,\\
g_{\ep_1,\dots,\ep_n}(\beta_1,\dots,\beta_n)&=&
\bigl(\delta_
{\sum_{j=1}^n \ep_j,-2}\bigr)
\prod_{\ell\in \L}\oint_{C_\ell} {d\al_\ell\over 2 \pi i}
\prod_{1\leq j<\ell \leq n}
(\al_\ell-\beta_j+{\pi i\over 2})\\
&&\times \prod_{1\leq \ell <j\leq n}
(\beta_j-\al_\ell+{\pi i\over 2})
\prod_{j,\ell}\Gamma(-{1\over 4}+
{\al_\ell-\beta_j\over 2 \pi i }) \Gamma(-{1\over 4}-
{\al_\ell-\beta_j\over 2 \pi i })\\
&&\times\prod_{\ell\in \L} \sinh(\al_\ell) 
\prod_{\ell<\ell^\prime}(\al_\ell-\al_{\ell^\prime}+\pi i)
\sinh(\al_\ell-
\al_{\ell^\prime}).
\ea
Moreover, $\Gamma(x)$   is the usual gamma
function, and 
each contour $C_\ell$ encloses only the poles $\alpha_\ell=
\beta_j-(2m+1)i\pi/2, m\geq 1$ and $\alpha_\ell=\beta_j+
i\pi/2$, with $1\leq j\leq n$, and is oriented anticlockwise.  
The set $\L$ is defined by
\be
\L=\{j,\quad s.t.\quad \ep_j=+1,\quad and \quad
\sum_{i=1}^n \ep_i=-2,\quad 1\leq j\leq n\}.
\ee 
Restricting
to the first Brillouin zone, that is,  
$0\leq k\leq 2\pi$, 
and performing the two delta integrations we find that 
$S^{i,+-}_n(w, k-\pi)$ is independent of $i$ 
and reads for $n\geq 2$
\bac
S^{i,+-}_n(w, k-\pi)&=&
C_n \int_{-\pi}^{0}dp_3\dots
\int_{-\pi}^0 dp_n\sum_{(\bar{p_1},\bar{p_2})} 
\Theta(W_u-W)\Theta(W-W_l)\\
&&\times {f(\beta(\bar p_1),\beta(\bar
p_2),\beta(p_3),\dots,\beta(p_n)) \sum_{\ep_1,\dots,\ep_n} 
|g_{\ep_1,\dots,\ep_n}(\beta(\bar p_1),\beta(\bar
p_2),\beta(p_3),\dots,\beta(p_n))|^2\over
\sqrt{W_u^2-W^2}}.
\label{exact}
\ea
Here, $\Theta$ is the Heaviside step function and
\bac
&&f(\beta(p_1),\beta(p_2),\beta(p_3),\dots,\beta(p_n))=
\prod_{1\leq j<j^\prime\leq
n}|A_{-}(\beta(p_{j^\prime})-
\beta(p_j))|^2,\\
&&C_n= 
{2^{(-3n^2+2n+4)/2} \pi^{7n(2-n)/4}  \over
n! \Gamma(1/4)^{n(n-2)}
|A_{-}(i\pi/2)|^{n(n-2)}},\\
&&W=w+\pi\sum_{j=3}^n \sin(p_j),\\
&&K=k+\sum_{j=3}^n p_j,\\
&&W_u=2 \pi |\sin(K/2)|,\\
&&W_l= \pi |\sin K|.
\ea
Moreover for fixed $W$ and $K$, the sum 
$\sum_{(\bar{p_1},\bar{p_2})}$ is over all  pairs  
$(\bar{p_1},\bar{p_2})$, solutions to the energy-momentum
conservation laws:
\bac
W&=&-\pi(\sin(\bar{p_1})+
\sin(\bar{p_2})),\\
K&=&-\bar{p_1}-\bar{p_2}.
\ea
The total $n$-spinon contribution is obtained through
\be
S^{+-}_n(w,k-\pi)=\sum_{i=0}^1 S^{i,+-}_n(w,k-\pi)=2 S^{i,+-}_n
(w, k-\pi),
\ee
from which we derive all the non-vanishing components of the
DSF as:
\be
S^{\mu\mu}_n(w,k-\pi)=2S^{+-}_n(w,k-\pi), \quad \mu=x,y,z.
\ee
The 2-spinon case, i.e.,  $n=2$, 
is particularly interesting 
because all contour integrals and  the
complicated $g$ function in (\ref{exact}) drop out.
Furthermore, because 
just the four pairs
$(\bar{p_1},\bar{p_2})$, $(\bar{p_2},\bar{p_1})$,
$(-\bar{p_1}-\pi,-\bar{p_2}-\pi)$ and 
$(-\bar{p_2}-\pi,-\bar{p_1}-\pi)$
contribute, and moreover equally, to the sum $\sum_{(\bar{p_1},\bar{p_2})}$,
$S^{+-}_2$ simplifies drastically to \cite{Boual96}
\bac
S^{+-}_2(w, k-\pi)&=&
8C_2  {\Theta(w_u-w)\Theta(w-w_l) 
|A_-(\beta(\bar{p_1})-\beta(\bar{p_2})|^2\over
\sqrt{w_u^2-w^2}},\\
&=&{\Theta(w_u-w)\Theta(w-w_l) 
|A_-(\beta(\bar{p_1})-\beta(\bar{p_2})|^2\over
\sqrt{w_u^2-w^2}}, 
\ea
where now 
\bac
w&=&-\pi(\sin(\bar{p_1})+
\sin(\bar{p_2})),\\
k&=&-\bar{p_1}-\bar{p_2},\\
w_u&=&2 \pi \sin(k/2),\\
w_l&=&\pi |\sin k|.
\ea
Let us mention that $S_2^{+-}(w,k-\pi)$ has been expressed in 
Ref. \cite{Karal96} just
in terms of the energy and momentum transfer through 
the
following relation: 
\be
\beta(\bar{p_1})-\beta(\bar{p_2})=2
cosh^{-1}\sqrt{{w_u^2-w_l^2\over w^2-w_l^2}}.
\ee

It is also interesting to note that $S^{+-}_2(w,k-\pi)$ vanishes 
at the upper boundary 
$w=w_u$ since $A_-(0)=0$. This is consistent with 
both facts that
at the latter boundary $\bar{p_1}=\bar{p_2}$ and that 
the spinons
obey the Pauli exclusion principle common to fermions. 
More detailed
analysis of the singular behaviour of $S^{+-}_2(w,k-\pi)$ 
at both 
boundaries 
$w=w_u$ and $w=w_l$ can be found in Ref. \cite{Karal96}.
Finally, let us mention that it would be particularly interesting
to evaluate the contour integrals that appear in the $g$
function of relation
(\ref{exact}), and especially for
$n=4$. The latter case is also useful  in optical experiments 
\cite{LoEd96}.

\vspace{0.25truein}

\section{Aknowlegements}

The work of A.H.B. is supported by the NSF Grant \# PHY9309888.
We are  grateful to Abada, Couture, Kacir, Karbach, Lorenzana, 
Saint-Aubin and Shrock for stimulating 
discussions and/or interesting comments.

\pagebreak


\begin{thebibliography}{10} 
 
\bibitem{Bet31}H. Bethe,
\newblock {Z. Phys.} {\bf 71}, 205 (1931).

\bibitem{Hul38} L. Hulth\'en,
\newblock {Arkiv Mat. Astron. Fysik A11} {\bf 26}, 1 (1938).

\bibitem{LiMa62} E.H. Lieb and D.C. Mattis,
\newblock{J. Math. Phys.} {\bf 3}, 749 (1962).

\bibitem{ClPe62} J. des Cloizeaux and J.J. Pearson,
\newblock {Phys. Rev.} {\bf 128}, 2131 (1962).

\bibitem{Gri64}R.B. Griffiths,
\newblock{Phys. Rev. } {\bf 133}, A768 (1964). 

\bibitem{YaYa66}C.N. Yang and C.P. Yang,
\newblock{Phys. Rev. } {\bf 150}, 321 (1966); {\bf 150}, 327
(1966); {\bf 151}, 258 (1966).

\bibitem{Gau71}M. Gaudin,

\newblock{Phys. Rev. Lett.} {\bf 26}, 1301 (1971).

\bibitem{Tak71}M. Takahashi,
\newblock{Prog. Theor. Phys.} {\bf 46}, 401 (1971).

\bibitem{FaTa79} L.A. Takhtajan and L.D. Faddeev,
\newblock { Russ. Math. Surveys} {\bf 34}, 11 (1979).

\bibitem{Babal83} O. Babelon, H.J. de Vega and C.M. Viallet,
\newblock {Nuc. Phys.} {\bf B220}, 13 (1983).


\bibitem{Mulal81} G. M\"uller, H. Thomas, H. Beck, and J.C.
Bonner,
\newblock { Phys. Rev.} {\bf B24}, 1429 (1981).

\bibitem{Mul82}G. M\"uller,
\newblock{Phys. Rev. } {\bf B26}, 1311 (1982). 

\bibitem{Rolal86} J.M.R. Roldan, B.M. McCoy and J.H.H. Perk,
\newblock { Physica } {\bf 136A}, 255 (1986).



\bibitem{Kor93}V.E. Korepin, A.G. Izergin, and N.M. Bogoliubov,
\newblock {\it Quantum Inverse Scattering 
Method and Correlation Functions},
\newblock Cambridge University Press, (1993).

\bibitem{Daval93}O. Davies, O. Foda, M. Jimbo, T. Miwa, and
A. Nakayashiki,
\newblock {\it Comm. Math. Phys.} 151 (1993) 89.




\bibitem{JiMi94}M. Jimbo and T. Miwa,
\newblock {\it Algebraic Analysis of Solvable Lattice Models},
\newblock {American Mathematical Society, (1994)}.

\bibitem{Fleal95}A. Fledderjohann, M. Karbach, K.-H. M\"utter and
P. Wielath,
\newblock{J. Phys.: Condens. Matter} {\bf 7}, 8993 (1995).

\bibitem{Fle96}A. Fledderjohann, M. Karbach and K.-H. M\"utter,
\newblock{Phys. Rev.} {\bf B53}, 11543 (1996).

\bibitem{Fleal96}A. Fledderjohann, C. Gerhardt, K.-H. M\"utter, 
A. Schmitt, and M. Karbach,
\newblock{Preprint cond-mat/9604085}.

\bibitem{Ten95} D.A. Tennant, R.A. Cowley, 
S.E. Nagler, and A. M. Tsvelik,
\newblock { Phys. Rev.} {\bf B52}, 13368 (1995).

\bibitem{Tenal95} D.A. Tennant, S.E. Nagler, D.Welz,  G.
Shirane, and K. Yamada,
\newblock { Phys. Rev.} {\bf B52}, 13381 (1995).

\bibitem{Boual96}A.H. Bougourzi, M. Couture and M. Kacir,
\newblock  Phys. Rev. {\bf B54}, 12669 (1996).

\bibitem{Karal96}M. Karbach, G. M\"uller, A.H. Bougourzi,
A. Fledderjohann and K.-H. M\"utter,
\newblock{Preprint ITP-SB-96-27}, Stony Brook,  1996.
To appear in Phys. Rev. B. 1997.

\bibitem{Fledd96} A. Fledderjohann, K.-H. M\"utter, M.
Karbach, and G. M\"uller
\newblock{Preprint WUB 96 cond-mat/9607102}.

\bibitem{Nie67} Th. Niemeijer, 
\newblock { Physica} {\bf 36}, 377 (1967).



\bibitem{LoEd96}J. Lorenzana and R. Eder,
\newblock {Preprint cond-mat/9608035.}
\end{thebibliography}
\end{document}